\documentclass{mcom-l}

\newcommand{\Z}{\mathbb{Z}}
\newcommand{\N}{\mathbb{N}}

\DeclareMathOperator{\ord}{ord}

\DeclareMathOperator{\integ}{int}

\usepackage{amssymb}
\usepackage{tikz}
\usetikzlibrary{arrows,backgrounds}
\usepackage{color}
\usepackage{algpseudocode}
\algtext*{EndWhile}
\algtext*{EndFor}
\algtext*{EndIf}

\copyrightinfo{}{}

\newtheorem{thm}{Theorem}[section]
\newtheorem{cor}[thm]{Corollary}
\newtheorem{lem}[thm]{Lemma}
\newtheorem{algorithm}[thm]{Algorithm}

\theoremstyle{definition}
\newtheorem{defin}[thm]{Definition}

\numberwithin{equation}{section}

\usepackage{geometry}
\begin{document}
\baselineskip=15pt
\title[Improving deterministic factorization]{A time-space tradeoff for\\ Lehman's deterministic integer factorization method}


\author[M. Hittmeir]{Markus Hittmeir}
\address{}
\curraddr{SBA Research, Floragasse 7, A-1040 Vienna}
\email{mhittmeir@sba-research.org}
\thanks{SBA Research (SBA-K1) is a COMET Centre within the framework of COMET – Competence Centers for Excellent Technologies Programme
	and funded by BMK, BMDW, and the federal state of Vienna. The COMET Programme is managed by FFG}


\date{}

\dedicatory{}

\begin{abstract}
Fermat's well-known factorization algorithm is based on finding a representation of natural numbers $N$ as the difference of squares. In 1895, Lawrence generalized this idea and applied it to multiples $kN$ of the original number. A systematic approach to choose suitable values for $k$ was introduced by Lehman in 1974, which resulted in the first deterministic factorization algorithm considerably faster than trial division. In this paper, we construct a time-space tradeoff for Lawrence's generalization and apply it together with Lehman's result to obtain a deterministic integer factorization algorithm with runtime complexity $O(N^{2/9+o(1)})$. This is the first exponential improvement since the establishment of the $O(N^{1/4+o(1)})$ bound in 1977.
\end{abstract}

\maketitle

\section{Introduction}
We consider the problem of computing the prime factorization of natural numbers $N$. There is a large variety of probabilistic and heuristic factorization methods achieving subexponential complexity. We refer the reader to the survey \cite{Len} and to the monographs \cite{Rie} and \cite{Wag}. The focus of the present paper is a more theoretical aspect of the integer factorization problem, which concerns \emph{deterministic} algorithms and their rigorous analysis. We will describe runtime complexities by using the bit-complexity model of the multitape Turing machine \cite{Pap}. 

In \cite{Str}, Strassen used fast polynomial multiplication and multipoint evaluation techniques to establish a deterministic and rigorous factorization algorithm running in time $O(N^{1/4+o(1)})$ bit operations. Since the publication of Strassen's approach in 1977, there have been a few refinements of  the runtime complexity. However, none of these improvements has been able to reduce the exponent $1/4$. The best of these bounds has been proved in \cite{Hit} and is given by
$\widetilde{O}\left(N^{1/4}\exp(-C\log N/\log\log N)\right)$ for a positive constant $C$, where the $\widetilde{O}$-notation is used to omit logarithmic factors. The contribution of this paper is to finally break the $1/4$-exponent threshold for deterministic integer factorization by proving the following theorem.

{\thm{There exists an algorithm to deterministically compute the prime factorization of any natural number $N$ in 
$
	O\left(\emph{\textsf{M}}_{\integ}\left(N^{2/9}\log^{5/3} N\right)\log N\right)
$
bit operations.}\label{t1}}
\vspace{12pt}

${\textsf{M}}_{\integ}(k)$ denotes the cost of multiplying two $\lceil k\rceil$-bit integers. Following \cite[p.1782]{BosGauSch}, we assume that
\begin{align}\label{e5}
\frac{\textsf{M}_{\integ}(k)}{k}\leq \frac{\textsf{M}_{\integ}(k')}{k'} \text { if $k\leq k'$ and }
\textsf{M}_{\integ}(kk')\leq k^2 \textsf{M}_{\integ}(k')
\end{align}
throughout the paper. Due to recent improvements (\cite{HarHoe}), ${\textsf{M}}_{\integ}(k)$ may be bounded by $O(k\log k)$. 

For the sake of simplicity, we will now consider $N=pq$, where $p$ and $q$ are distinct primes. In order to prove Theorem \ref{t1}, we will establish a time-space tradeoff for Fermat's well-known factorization method and its generalizations by Lawrence (\cite{Law}) and Lehman (\cite{Leh}). Fermat's method is a special-purpose algorithm first described in 1643. The idea of the procedure is based on the fact that $(p+q)^2-(q-p)^2=4N$. If the prime difference $\Delta:=q-p$ is small, the number $S:=p+q$ is close to $2\sqrt{N}$. In fact, one proves that
\[
	0<S-2\sqrt{N}<\frac{\Delta^2}{4N^{1/2}}. \label{e1}
\]
In the algorithm, we hence try to find $S$ by computing $(z+\lceil 2\sqrt{N}\rceil)^2-4N$ for $z=0,1,2,3\ldots$ and checking if one of those numbers is equal to a square. If this procedure indeed yields an equality of the form $(z+\lceil 2\sqrt{N}\rceil)^2-4N=y^2$, we may obtain a proper factor of $N$ by computing $\gcd(z+\lceil 2\sqrt{N}\rceil-y,N)$. Lawrence's generalization works in a similar manner, but is based on finding linear combinations $ap+bq$ of the prime factors of $N$. In 1974, Lehman published a deterministic factorization method running in time $\widetilde{O}(N^{1/3})$, thus being significantly faster than trial division. The algorithm applies Lawrence's approach in a systematic manner and relies on suitable lower and upper bounds for the linear combinations $ap+bq$ (see Theorem \ref{t2}). All pairs $(a,b)$ in a certain range are checked,  and  the bounds for the corresponding linear combinations vary depending on the size of the value of $ab$. 

In our improvement, we utilize Lehman's theorem and elaborate on an idea from \cite{Hit}. There, it was shown that $\alpha^{S}\equiv \alpha^{N+1}\mod N$ holds for every $\alpha$ coprime to $N$. Subsequently, this congruence has been used to search for $S=p+q$. In Lemma \ref{l1} of the present paper, we will extend this statement to the linear combinations $ap+bq$ used by Lehman. As a consequence, we are able to construct two sets $\mathcal{B}$ and $\mathcal{S}$ which are disjoint modulo $N$, but not disjoint modulo $p$ or $q$. In order to find the corresponding collision and, hence, a proper factor of $N$, we will adapt the already mentioned factorization scheme by Strassen. The application of efficient polynomial arithmetic techniques together with Lehman's bounds lead to the runtime complexity stated in Theorem \ref{t1}.

The remainder of the paper is organized as follows: In Section \ref{s1}, we discuss two theorems important to our improvement, one of which is Lehman's result. In Section \ref{s2}, we consider our adaptation of Strassen's factorization scheme and introduce the concept of revealing subsets. Section \ref{s3} contains the core idea of this paper, namely the time-space tradeoff for Lawrence's algorithm. In Section \ref{s4}, we explain our strategy for applying Lehman's bounds. Finally, Section \ref{s5} puts all the pieces together and finishes the proof.

\section{Preliminary Work}\label{s1}
This short section discusses two preliminary results our improvement is based on. Throughout the paper, we will use the notations $\N=\{1,2,3,\ldots\}$ and $\Z_n:=\Z/n\Z$ for $n\in\N$.
The first and most important ingredient of our algorithm is the following theorem, which was proved in \cite{Leh}. Lehman applied it with $\eta\in O(N^{1/3})$ to obtain a new factorization technique running in time $\widetilde{O}(N^{1/3})$. Compared to trial division, this was a major improvement.

{\thm{\emph{[Lehman, 1974]} Suppose that $N\in\N$ is odd and $\eta$ is an integer such that $1\leq \eta<N^{1/2}$. If $N=pq$, where $p$ and $q$ are primes and 
		\[
		\left(N/(\eta+1)\right)^{1/2}< p \leq N^{1/2},
		\]
		then there are non-negative integers $x$, $y$ and $k$ such that
		\begin{align}
			&x^2-y^2=4kN,\text{ } 1\leq k\leq \eta, \notag\\
			&0\leq x-(4kN)^{1/2}\leq \frac{N^{1/2}}{4k^{1/2}(\eta+1)}. \label{e3}
		\end{align}
		These integers are of the form $x=ap+bq$, $y=|ap-bq|$ and $k=ab$, where $a,b\in\N$.
}\label{t2}}
\vspace{12pt}

Lehman used Lawrence's approach systematically and searched for linear combinations of the prime divisors of $N$. The idea is to divide the interval $[0,1]$ into parts, where each part corresponds to a fraction $a/b$. Considering many such fractions, the goal is to find the best approximation for $q/p$. This implies that the value of $|ap-bq|$ is sufficiently small such that the inequality $(\ref{e3})$ holds true. The shape of $x,y$ and $k$  is shown in (4.6) of Lehman's proof in \cite[p.641]{Leh} and is not part of the original formulation of the theorem. However, it is crucial for our application of the result.

Our second ingredient is a deterministic algorithm for finding elements of large order modulo $N$. We will need such an element in our application of the already mentioned Lemma \ref{l1}. The following result has been proved in \cite[Theorem 6.3]{Hit}.

{\thm{\emph{[Hittmeir, 2018]} Let $N\in\N$ and $\delta$ be an integer such that $N^{2/5}\leq \delta\leq N$. There exists an algorithm that either returns some $\alpha\in\Z_N^*$ with $\ord_N(\alpha)>\delta$, or some nontrivial factor of $N$, or proves $N$ prime. Its runtime complexity is bounded by 
		\[
		O\left(\frac{\delta^{1/2}\log^2 N}{\sqrt{\log\log \delta}} \right)
		\]
		bit operations.
}\label{t3}}
\vspace{12pt}

In this algorithm, we apply the standard babystep-giantstep procedure to compute small orders of elements $\alpha\in\Z_N^*$. If $m:=\ord_N(\alpha)$ is actually found, we try to find a factor of $N$ via $\gcd(N,\alpha^{m/r}-1)$ for every $r\mid m$. If this fails, we know that $m\mid p-1$ for every prime divisor $p$ of $N$. Repeating this process for various values of $\alpha$, we either find an element of sufficiently large order or obtain enough information about the factorization of $p-1$ for $p\mid N$ such that factoring $N$ directly is feasible.

\section{Revealing Subsets} \label{s2}
We briefly recall Strassen's idea for factoring natural numbers. Let $N\in\N$ and set $d:=\lceil N^{1/4}\rceil$. We want to compute subproducts of the product $\lceil N^{1/2}\rceil !$ to find a factor of $N$. For this task, the polynomial 
$
f=(X+1)(X+2)\cdots (X+d)
$
is computed modulo $N$ and evaluated at the points $0,d,2d,\ldots,(d-1)d$ by using fast polynomial arithmetic techniques. In the next step, $g_i:=\gcd(f(id),N)$ is computed for $i=0,\ldots,d-1$. If $N$ is not a prime, one of these GCDs is not equal to $1$. If $g_{i_0}=N$, we obtain a nontrivial factor of $N$ by computing $\gcd(i_0d+j,N)$ for $j=1,\ldots,d$.

In this paper, we will apply Strassen's idea in a generalized setting. We will still use linear polynomials as factors of $f$, but with different sets of zeros. Similarly, we will use different sets of evaluation points for computing the GCDs. The following definition clarifies the conditions under which such sets will reveal a factor of $N$ in the same way the original procedure of Strassen does. Subsequently, we formulate the corresponding algorithm.

{\defin{Let $N\in\N$. A pair of subsets $\mathcal{B}$ and $\mathcal{S}$ of $\Z_N^*$ is called \emph{revealing} if the following two conditions hold:
		\begin{enumerate}
			\item{$\mathcal{B}$ and $\mathcal{S}$ are disjoint modulo $N$.}
			\item{If $N$ is composite, then there is a nontrivial divisor $v$ of $N$ such that $\mathcal{B}$ and $\mathcal{S}$ are not disjoint modulo $v$.} 
		\end{enumerate}
}\label{d1}}

\begin{algorithm}
	\emph{Input:} A natural number $N$ and a pair of disjoint subsets $\mathcal{B}$ and $\mathcal{S}$ of $\Z_N^*$. We denote $\mathcal{B}=\{b_i:1\leq i\leq \beta\}$ and $\mathcal{S}=\{s_j:1\leq j\leq \sigma\}$.
	
	\begin{algorithmic}[1]
		\State Compute
		$
		f:=\prod_{i=1}^\beta (X-b_i). 
		$
		\State Compute $f(s_1),f(s_2),\ldots,f(s_\sigma)$.
		\For{$j=1,\ldots, \sigma$}
		\State Compute $\gamma_j:=\gcd(f(s_j),N)$.
		\If {$1<\gamma_{j}<N$} return $\gamma_{j}$.
		\EndIf
		\If{$\gamma_{j}=N$}
		\For{$i=1,\ldots,\beta$}
		\State  Compute $\gamma_{j,i}:=\gcd(s_j-b_i,N)$.
		\If {$\gamma_{j,i}>1$} return $\gamma_{j,i}$.
		\EndIf
		\EndFor
		\EndIf
		\EndFor
		\State Return ``No factor found''
	\end{algorithmic}
	\label{a1}
\end{algorithm}

In the subsequent sections, we will employ Algorithm \ref{a1} with different choices for $\mathcal{B}$ and $\mathcal{S}$. We now prove correctness and analyze the runtime of the procedure.

{\lem{Algorithm \ref{a1} runs in time
		$
		O\left(\emph{\textsf{M}}_{\integ}\left(\delta\log N\right)\log N\right),
		$
		where $\delta:=\max\{\beta,\sigma\}$. If $N$ is composite and $\mathcal{B}$ and $\mathcal{S}$ are revealing, it returns a nontrivial factor of $N$.
}\label{l3}}
\begin{proof}
	We first prove correctness. Let $N$ be composite and $\mathcal{B}$ and $\mathcal{S}$ be revealing subsets. There are elements $b_{i'}\in \mathcal{B}$ and $s_{j'}\in\mathcal{S}$ such that $b_{i'}\equiv s_{j'}\mod v$ for some nontrivial divisor $v$ of $N$. As a result, one easily observes that $\gamma_{j'}>1$ in Step 4 of the algorithm. If $\gamma_{j'}<N$, the algorithm returns a correct output in Step 5. If $\gamma_{j'}=N$, then $\gamma_{j',i'}$ in Step 8 will be a nontrivial factor of $N$, since $\mathcal{B}$ and $\mathcal{S}$  are disjoint modulo $N$. In this case, the algorithm terminates in Step 9. 
	
	We now discuss the runtime complexity. Since $\mathcal{B}$ and $\mathcal{S}$ are subsets of $\Z_N^*$, we have $\delta<N$. In the Steps 1 and 2, we use the subproduct-tree algorithm to compute $f$ and apply multipoint-evaluation techniques to obtain the the values $f(s_j)$ for $j=1,\ldots,\sigma$. Detailed information about these procedures can be found in \cite[Chap.10]{GerGat}. The complexity is bounded by $O(\textsf{M}(\delta)\log N)$, where $\textsf{M}(k)$ denotes the cost of multiplying two polynomials of degree less than $k$ in $\Z_N^*$. Applying Lemma 3 in \cite{BosGauSch} with $m=1$, this reduces to a bit-complexity of $O\left(\textsf{M}_{\integ}\left(\delta\log N\right)\log N\right)$. We now consider the complexity of the GCD computations in the Steps 3-10. The fact that $\mathcal{B}$ and $\mathcal{S}$  are disjoint modulo $N$ implies that we only reach the loop in Step 7 in cases where the algorithm terminates by returning a factor in Step 9. As a result, the number of computed GCDs is bounded by $O(\delta)$. Since the (extended) GCD of two bit-size $k$ integers can be obtained in time $O(\textsf{M}_{\integ}(k)\log k)$ (see \cite[p.1783]{BosGauSch}), the complexity of the Steps 3-10 is negligible.
\end{proof}

One observes that we obtain the original approach of Strassen as a special case. We will apply it in form of the following corollary.

{\cor{Let $N\in\N$ and $\Delta\in\N$ with $\Delta\leq N^{1/2}$. We can find a nontrivial divisor $\ell$ of $N$ such that $\ell\leq \Delta$, or prove that no such divisor exists, in
		$
	O\left(\emph{\textsf{M}}_{\integ}\left(\Delta^{1/2}\log N\right)\log N\right)
		$
		bit operations.
}\label{c1}}

\begin{proof}
Assume that $N$ is sufficiently large and let $d=\lceil \Delta^{1/2}\rceil$. We note that $1\leq d<N$ and the subsets $\mathcal{B}=\{-i: 0\leq i\leq d-1\}$ and $\mathcal{S}=\{dj:1\leq j\leq d\}$ of $\Z_N^*$ are disjoint. One easily shows that they are revealing if and only if there is a nontrivial divisor $\ell$ of $N$ with $\ell\leq \Delta$. We now apply Lemma \ref{l3}.	
\end{proof}

Prior to the application of Algorithm \ref{a1}, we always have to make sure that the input subsets are revealing in the sense of Definition \ref{d1}. In order to prepare two lists of elements in $\Z_N^*$ for the application of Algorithm \ref{a1}, we will use the following lemma. Note that the notion of indexed lists refers to lists where each element is stored together with an index, which we assume to be an integer or a tupel of integers. 

{\lem{Let $N\in\N$ and $\mathcal{L}_1$ and $\mathcal{L}_2$ be two indexed lists of elements in $\Z_N^*$ and indices of bit size $O(\log N)$. Moreover, assume that the elements of $\mathcal{L}_1$ are all distinct. We may compute the set $\mathcal{M}$ of all the pairs of indices $(i,j)$ for which the element with index $i$ in $\mathcal{L}_1$ is equal to the element with index $j$ in $\mathcal{L}_2$. The bit-complexity is bounded by $O(\delta \log^2 N)$, where $\delta$ is the maximum of the lengths of the lists.
}\label{l4}}

\begin{proof}
	In the multitape Turing machine model, we use merge sort to solve this task. We denote every number in the lists as string of at most $O(\log N)$ bits and sort both $\mathcal{L}_1$ and $\mathcal{L}_2$ by performing at most $O(\delta\log \delta)$ comparisons, yielding a bit-complexity bounded by $O(\delta\log^2 N)$. After applying merge sort, it is easy to see that one is able to find all matches in the lists in $O(\delta\log N)$ bit operations. For more detailed information on merge sort, we refer the reader to \cite[Chap. 2.2]{SedWay}. For the claimed complexity, see Proposition F in \cite[Page 272]{SedWay}.
\end{proof}

\section{Time-Space Tradeoff for Lawrence's algorithm \label{s3}}
Let $N\in\N$ and $u$ and $v=N/u$ be nontrivial divisors of $N$. As discussed in the introduction, Fermat's factorization algorithm searches for $S=u+v$ by checking if $i^2-4N$ is a square for each $i\geq\lceil2\sqrt{N}\rceil$. Lawrence (\cite{Law}) first noted that, since the equality $(au+bv)^2-(au-bv)^2=4abN$ holds true for all $a,b\in\N$, we can try to find any linear combination of $u$ and $v$, not only $S$. For each $i\geq\lceil2\sqrt{abN}\rceil$, we check if $i^2-4abN$ is a square number. As in Fermat's original approach, we have
\begin{equation}
0<au+bv-2\sqrt{abN}<\frac{(au-bv)^2}{4\sqrt{abN}}, \label{e1}
\end{equation}
hence finding $au+bv$ works best if $|au-bv|$ is relatively small or, in other words, if $a/b$ is a good approximation of $v/u$. The following lemma shows how to factorize $N$, given that we know $au+bv$.

{\lem{Let $N\in\N$ and $u$ and $v=N/u$ be unknown nontrivial divisors of $N$. Given as input positive numbers $a,b,L\in O(N)$ such that $a$ and $b$ are coprime to $N$, we may test if $L=au+bv$ or $L=bu+av$, and if so compute $u$ and $v$, in $O(\emph{\textsf{M}}_{\integ}(\log N))$ bit operations.}\label{l2}}

\begin{proof}
We compute the roots $r_1$ and $r_2$ of the polynomial
$
X^2-LX+abN.
$	
One easily checks that either $r_1/a$ or $r_1/b$ is equal to one of the factors of $N$. Since all involved quantities have $O(\log N)$ bits, we achieve the claimed runtime complexity by using the quadratic formula.
\end{proof}

In this section, we introduce a time-space tradeoff for Lawrence's extension for semiprime numbers $N$. The extension is based on the following lemma, which may be considered as a generalization of the fact that $\alpha^{p+q}\equiv \alpha^{N+1}\mod N$ for every $\alpha\in\Z_N^*$.
{\lem{
Let $N$ be semiprime with the distinct factors $p$ and $q$. Furthermore, let $a,b\in\Z$, $\alpha\in\Z_N^*$ and set $t:= \alpha^{bN+a}\pmod N$. We have $\alpha^{ap+bq}\equiv t \mod p$ and $\alpha^{bp+aq}\equiv t \mod q$.
}\label{l1}}

\begin{proof}
The first congruence follows from the fact that $p\equiv 1 \mod p-1$ and, hence, 
\[
ap+bq\equiv a+bpq \equiv a+bN \mod p-1.
\]
The proof of the second congruence is similar.	
\end{proof}

Setting $t:= \alpha^{bN+a}\pmod N$ introduces an asymmetry of $a$ and $b$ and defines which congruence in Lemma \ref{l1} holds modulo which prime factor. The following arguments show how to find $ap+bq$ by considering the first congruence modulo $p$. However, they may easily be modified to fit the search for $bp+aq$ via the congruence modulo $q$. In the subsequent theorem, we search for both quantities at once.

Let $x_0=ap+bq-\lceil2\sqrt{abN}\rceil$ and $t_0:=\alpha^{bN+a-\lceil2\sqrt{abN}\rceil}\pmod N$. From Lemma \ref{l1}, it follows that  
$
\alpha^{x_0}\equiv t_0 \mod p.
$
Now let $\Lambda\in \N$ such that $x_0<\Lambda$. For example, $\Lambda$ may be taken from inequality (\ref{e1}) or Lehman's bound in Theorem \ref{t2}. Let $m=\lceil \Lambda^{1/2}\rceil$ and write $x_0=ml_0+r_0$ for unknown numbers $l_0,r_0\in \{0,1,\ldots,m-1\}$. We then have
\begin{equation}
\alpha^{r_0}\equiv \alpha^{-ml_0}t_0 \mod p.\label{e2}
\end{equation}
Our goal is to solve (\ref{e2}) in order to find $l_0$ and $r_0$ and obtain a proper factor of $N$. However, (\ref{e2}) is a congruence modulo an unknown prime factor. We will hence employ the approach discussed in Section \ref{s2}, which allows us to prove the following result.

{\thm{Let $N$ and $\Lambda$ be natural numbers. There is an algorithm that takes as input $\alpha\in\Z_N^*$ such that $\ord_N(\alpha)>\lceil\Lambda^{1/2}\rceil$ and positive integers $a,b\in O(N)$ which are coprime to $N$. Its runtime complexity is bounded by
\[
O\left(\emph{\textsf{M}}_{\integ}\left(\Lambda^{1/2}\log N\right)\log N\right).
\]
If $N$ is semiprime with distinct factors $p$ and $q$, and $ap+bq-\lceil2\sqrt{abN}\rceil< \Lambda$ or $bp+aq-\lceil2\sqrt{abN}\rceil<\Lambda$ holds, the algorithm returns $p$ and $q$. }\label{t4}}

\begin{proof}
	Put $m=\lceil\Lambda^{1/2}\rceil$ and let $x_0=ap+bq-\lceil2\sqrt{abN}\rceil$ and $x_1=bp+aq-\lceil 2\sqrt{abN}\rceil$. Our assumptions imply that $x_0< \Lambda$ or $x_1 <\Lambda$ holds, and we will use Lemma \ref{l1} to search for them simultaneously. W.l.o.g., we assume that $x_0< \Lambda$ and $x_0=ml_0+r_0$ for unknown numbers $l_0,r_0\in \{0,1,\ldots,m-1\}$. Before applying Algorithm \ref{a1}, we perform D. Shanks' babystep-giantstep  method to search for matches of the form (\ref{e2}), but modulo $N$ instead of modulo $p$. We start by computing $t_0$ in time $O(\textsf{M}_{\integ}(\log N)\log N)$ via the Square-and-Multiply algorithm. Next, we compute the lists of the babysteps $\alpha^i \pmod N$ for $0\leq i \leq m-1$ and the giantsteps $\alpha^{-mj}t_0 \pmod N$ for $0\leq j \leq m-1$. Note that these elements should be stored together with their indices $i$ and $j$. It is easy to see that this can be done by performing $O(m)$ multiplications modulo $N$, which is asymptotically negligible. 
	
	In the next step of the algorithm, we compute 
	$
	g_i:=\gcd(N,\alpha^i -1\pmod N)
	$
	for $i=1,\ldots,m-1$. The computation of $m$ GCDs of numbers bounded by $O(N)$ can be done in $O(m\cdot \textsf{M}_{\integ}(\log N) \log \log N)$ bit operations, which is also negligible. Since we assumed $\ord_N(\alpha)>m$, it follows that $g_i<N$ for every $i$. As a result, we either find a factor of $N$ or we obtain $g_i=1$ for every $i$. If we do not find a factor, we know that $\ord_p(\alpha)\geq m$ and $\ord_q(\alpha)\geq m$. We now suppose this is the case. Since our assumptions on the order of $\alpha$ modulo $N$ imply that all the babysteps are distinct, we may apply Lemma \ref{l4} to compute the set $\mathcal{M}$ of pairs of indices corresponding to the matches between the babysteps and the giantsteps. The runtime is bounded by $O(m\log^2 N)$, which again is negligible. Moreover, we deduce that $|\mathcal{M}|\leq m$. For each solution $(i',j')$ in $\mathcal{M}$, we check if $mj'+i'=x_0$. This can be done by using 
	$
	mj'+i'+\lceil2\sqrt{abN}\rceil
	$
	as candidate for $L$ in Lemma \ref{l2}. Assume that we do not find $p$ and $q$. The lower bounds on the order of $\alpha$ modulo $p$ and modulo $q$ imply that $\alpha^{i'}$ is the only babystep that matches  $\alpha^{-mj'}t_0$ modulo the prime factors of $N$. It follows that $j'\neq l_0$. As a result, we delete the element $\alpha^{-mj'}t_0 \pmod N$ from the list of giantsteps. The overall runtime complexity for these applications of Lemma \ref{l2} can be bounded by $O(m\cdot \textsf{M}_{\integ}(\log N))$, and hence is negligible. 
	
	If we have not found a factor at this point, we define the set $\mathcal{B}$ consisting of the babysteps and the set $\mathcal{S}$ consisting of the remaining giantsteps. Note that the list of giantsteps may have contained several instances of the same element in $\Z_N^*$, stored together with different indices. In the sets $\mathcal{B}$ and  $\mathcal{S}$, we ignore multiplicities and do not longer keep track of the indices. Clearly, $\mathcal{B}$ and $\mathcal{S}$ are disjoint modulo $N$. Moreover, we point out that the set of giantsteps is not empty, as it must still contain the element $\alpha^{-ml_0}t_0 \pmod N$. It follows that the two sets are not disjoint modulo $p$ or not disjoint modulo $q$. We conclude that $\mathcal{B}$ and  $\mathcal{S}$ are a pair of revealing subsets of $\Z_N^*$. We hence may apply Lemma \ref{l3} and the claim follows.
\end{proof}

\section{Utilizing Lehman's bound} \label{s4}
We now discuss the application of Lehman's bound (\ref{e3}) in Theorem \ref{t2}, which is at the core of our improvement. In the following, we suppose that $N$ is a prime or a semiprime number. A procedure for reducing the factorization of any natural number to the factorization of primes and semiprimes will be discussed at the beginning of the next section. Assuming that $N=pq$, our goal is to find the number $x$ which satisfies the bound (\ref{e3}). As discussed in Section \ref{s1}, this number is of the shape $x=ap+bq$. We consider all possible values of $k=ab$ in the interval $[1,\eta]$ and apply the time-space tradeoff established in Theorem \ref{t4} to compute the factors of $N$. We have to account for the fact that said approach cannot search for all decompositions of a fixed $1\leq k \leq \eta$ at once. Instead, we have to go through all $(a,b)$ with $1\leq ab\leq \eta$. However, we will see that the effect on the runtime complexity is marginal. Moreover, we may assume that $a\leq b$, since our procedure still searches for $ap+bq$ and $bp+aq$ at once.

In addition, we split the interval $[1,\eta]$ into $[1,\xi]$ and $(\xi,\eta]$, employing different strategies for the pairs $(a,b)$ depending on whether $ab$ is in the first or the second range. This idea is based on the fact that Lehman's bound (\ref{e3}) is large for smaller values of $k\leq \xi$, which allows for an effective application of the procedure in Theorem \ref{t4} for every single pair $(a,b)$. However, for larger values of $k>\xi$, the bound (\ref{e3}) is decreasing rapidly. In this case, we are able to denote \emph{all remaining candidates} for $x$ as $\delta+\lceil 2\sqrt{kN}\rceil$ for some sufficiently small $\delta$, and we will require only one more run of an approach similar to the one taken in Theorem \ref{t4} to find $x$ and, hence, the prime factors $p$ and $q$.  
 
Both parameters $\xi$ and $\eta$ have a strong impact on the final runtime complexity and will be optimized later. Based on the explanations above, we now suggest the following algorithm for a time-space tradeoff of Lehman's algorithm.


\begin{algorithm}
	\emph{Input:} A natural number $N$ and two integers $\xi$ and $\eta$ such that $1\leq \xi\leq \eta<N^{1/4}$ and $\left(N/(\eta+1)\right)^{1/2}< p \leq N^{1/2}$. Moreover, an element $\alpha\in\Z_N^*$ such that $\ord_N(\alpha)>\lceil N^{1/2}/(\xi^{1/2}\eta)\rceil$. 
	
	\begin{algorithmic}[1]
		\State For every $(a,b)$ with $a\leq b$ and $ab\leq \xi$, apply Theorem \ref{t4} with Lehman's bound (\ref{e3}), i.e.
		\[
		\Lambda_{a,b}:=\left\lceil\frac{N^{1/2}}{(ab)^{1/2}\eta}\right\rceil> \frac{N^{1/2}}{4(ab)^{1/2}(\eta+1)}.
		\]
		If a factor is found, stop.
		\State Construct two indexed lists: $\mathcal{L}_1$ consists of $\alpha^i\pmod N$ with the indices $i=0,1,...,\lceil N^{1/2}/(\xi^{1/2}\eta))\rceil$. $\mathcal{L}_2$ consists of $\alpha^{bN+a-\lceil2\sqrt{abN}\rceil}\pmod N$ with the indices $(a,b)$, where we consider all $(a,b)$ satisfying $a\leq b$ and $\xi<ab\leq \eta$. 
		\State Compute $\gcd(N,\alpha^i-1\pmod N)$ for $i=0,1,...,\lceil N^{1/2}/(\xi^{1/2}\eta))\rceil$. If a factor is found, stop.
		\State Apply Lemma \ref{l4} to $\mathcal{L}_1$ and $\mathcal{L}_2$. For every element $(i',(a',b'))$ in the resulting set $\mathcal{M}$, apply Lemma \ref{l2} with 
		$
		L=i'+\lceil2\sqrt{a'b'N}\rceil
		$
		as candidate. If a factor is found, stop. If not, delete the element with index $(a',b')$ from $\mathcal{L}_2$.
		\State Let $\mathcal{B}$ be the set of the elements in $\mathcal{L}_1$ and $\mathcal{S}$ the set of the remaining elements in $\mathcal{L}_2$. Apply Algorithm \ref{a1} to $\mathcal{B}$ and $\mathcal{S}$.
	\end{algorithmic}
	\label{a2}
\end{algorithm}

{\lem{The bit-complexity of Algorithm \ref{a2} is bounded by 
\[
O\left(\emph{\textsf{M}}_{\integ}\left(\frac{N^{1/4}\xi^{3/4}\log \xi}{\eta^{1/2}}\log N\right)\log N+\emph{\textsf{M}}_{\integ}\left(\frac{N^{1/2}}{\xi^{1/2}\eta}\log N\right)\log N+\emph{\textsf{M}}_{\integ}\left(\eta\log\eta\log N\right)\log N\right).
\]
If $N$ is semiprime with distinct factors $p$ and $q$, it returns a prime factor of $N$. 
}\label{l5}}

\begin{proof}
In order to prove correctness, we let $\bar{a}$ and $\bar{b}$ such that $x:=\bar{a}p+\bar{b}q$ and $k:=\bar{a}\bar{b}$ in Theorem \ref{t2}. If $k\leq \xi$, one easily checks that $p$ and $q$ will be found in the corresponding run of Step 1. Just note that our assumptions imply $\ord_N(\alpha)>\lceil N^{1/2}/(\xi^{1/2}\eta)\rceil>\lceil N^{1/4}/\eta^{1/2}\rceil\geq\Lambda_{a,b}^{1/2}$ for all $(a,b)$, hence Theorem \ref{t4} is applied correctly. We also point out that it suffices to consider the pairs $(a,b)$ with $a\leq b$, since Theorem \ref{t4} searches for both $ap+bq$ and $bp+aq$ at once. 

We now assume that $\xi<k\leq \eta$. In this case, note that Lehman's bound (\ref{e3}) implies that we may write $x=\delta + \lceil 2\sqrt{kN} \rceil$ for some
\[
0\leq \delta \leq \frac{N^{1/2}}{4k^{1/2}(\eta+1)}<\left\lceil \frac{N^{1/2}}{\xi^{1/2}\eta}\right\rceil.
\]
As a result, we derive
$
\alpha^\delta\equiv \alpha^{\bar{b}N+\bar{a}-\lceil2\sqrt{kN}\rceil} \mod p.
$
In Step 2, we construct indexed lists modulo $N$ such that there is a collision modulo $p$ which corresponds to this congruence. The remainder of the algorithm and the following arguments are similar to the approach we took in the proof of Theorem \ref{t4}. The goal is to end up with a pair of revealing subsets in $\Z_N^*$. In Step 3, we either find a prime factor of $N$ or make sure that both $\ord_p(\alpha)$ and $\ord_q(\alpha)$ are greater or equal to $\lceil N^{1/2}/(\xi^{1/2}\eta)\rceil$. In Step 4, we search for matches between $\mathcal{L}_1$ and $\mathcal{L}_2$ modulo $N$ and check if the corresponding indices yield candidates that are equal to $x$. The lower bounds on the order of $\alpha$ modulo $p$ and modulo $q$ imply that at most one babystep matches each giantstep modulo these prime factors. Since we delete all matching elements from $\mathcal{L}_2$ in Step 4, the sets $\mathcal{B}$ and $\mathcal{S}$ in Step 5 are disjoint modulo $N$. Moreover, the fact that no factor has been found at this point implies that the element corresponding to $(\bar{a},\bar{b})$ is still in $\mathcal{S}$. Of course, the element $\alpha^\delta\pmod N$ is in $\mathcal{B}$, and the desired result follows.

We proceed by considering the runtime complexity of the algorithm. For Step 1, we consider Theorem \ref{t4} and the sum over all $(a,b)$ with $a\leq b$ and $ab\leq \xi$. Hence, we write
\begin{align*}
\sum_{a=1}^{\lfloor \xi^{1/2}\rfloor}\sum_{b=a}^{\lfloor \xi/a\rfloor}\textsf{M}_{\integ}\left(\Lambda_{a,b}^{1/2}\log N\right)\log N
&=\sum_{a=1}^{\lfloor \xi^{1/2}\rfloor}\sum_{b=a}^{\lfloor \xi/a\rfloor}\textsf{M}_{\integ}\left(\frac{N^{1/4}}{(ab)^{1/4}\eta^{1/2}}\log N\right)\log N\\
&\leq
\textsf{M}_{\integ}\left(\frac{N^{1/4}\log N}{\eta^{1/2}}\cdot \sum_{a=1}^{\lfloor \xi^{1/2}\rfloor}\frac{1}{a^{1/4}}\sum_{b=a}^{\lfloor \xi/a\rfloor}\frac{1}{b^{1/4}}\right)\log N,
\end{align*}
where we have used that $\textsf{M}_{\integ}(k)+\textsf{M}_{\integ}(k')\leq \textsf{M}_{\integ}(k+k')$ holds for all $k,k'$, which follows from the first assumption in (\ref{e5}). Now by considering 
\begin{align*}
\sum_{a=1}^{\lfloor\xi^{1/2}\rfloor}\frac{1}{a^{1/4}}\sum_{b=a}^{\lfloor \xi/a\rfloor}\frac{1}{b^{1/4}}&\leq \sum_{a=1}^{\lfloor \xi^{1/2}\rfloor}\frac{1}{a^{1/4}}\int_{0}^{\lfloor \xi/a\rfloor}t^{-1/4}dt\\ 
&\in O\left(\sum_{a=1}^{\lfloor \xi^{1/2}\rfloor}\frac{\xi^{3/4}}{a}\right)\subseteq O(\xi^{3/4} \log \xi),
\end{align*}
we are able to obtain the first summand in our claimed runtime complexity bound. For estimating the runtime complexity of the Steps 2-4, let $\kappa$ denote the maximum of the lengths of the two lists $\mathcal{L}_1$ and $\mathcal{L}_2$. We first note that the length of $\mathcal{L}_2$ is certainly no more than the number of pairs $(a,b)$ with $1\leq ab\leq \eta$, which is bounded by
\begin{align*}
\sum_{a=1}^{\lfloor \eta^{1/2}\rfloor}\sum_{b=1}^{\lfloor \eta/a\rfloor}1&\leq \sum_{a=1}^{\lfloor \eta^{1/2}\rfloor}\frac{\eta}{a}\in  O(\eta\log\eta).
\end{align*}
Therefore, it follows that
\[
\kappa\in O \left(\max\left\{\frac{N^{1/2}}{\xi^{1/2}\eta},\eta \log \eta\right\}\right).
\]
Since all exponents in Step 2 are in $O(N^2)$, the construction of the two lists via the well-known Square-and-Multiply algorithm is finished after at most $O (\kappa\cdot \textsf{M}_{\integ}\left(\log N\right)\log N)$ bit operations, which is negligible. The computation of the GCDs in Step 3 is finished after at most $O(\kappa\cdot \textsf{M}_{\integ}(\log N) \log \log N)$ bit operations. The complexity for applying Lemma \ref{l4} in Step 4 is $O(\kappa \log^2 N)$, which is also negligible due to the fact that $k\leq \textsf{M}_{\integ}(k)$ is true for all $k$. The same holds for the cost of applying Lemma \ref{l2} with all possible candidates, which may be bounded by $O(\kappa\cdot \textsf{M}_{\integ}(\log N))$. Finally, the runtime complexity for Step 5 is a result of Lemma \ref{l3}, which finishes the proof.	
\end{proof}

\section{Proof of Theorem \ref{t1}} \label{s5}
Let $N$ be any natural number. We start by performing the following preparation step: We apply Corollary \ref{c1} to $N$, putting $\Delta_0=\lceil N^{1/3}\rceil$. If any divisor of $N$ is found, we remove it and denote the resulting number by $N_1$. We then apply Corollary \ref{c1} to $N_1$, setting $\Delta_1=\lceil N_1^{1/3}\rceil$. We proceed in this way until no more divisors are found. The number of factors of $N$ obtained in this manner is bounded by $O(\log N)$. Since $\Delta_i\leq \lceil N^{1/3}\rceil$ for every $i\in\N_0$, the bit-complexity of this procedure is at most
\[
O\left(\textsf{M}_{\integ}\left(N^{1/6}\log N\right)\log^2 N\right).
\]
Compared to the complexity bound stated in Theorem \ref{t1}, this is asymptotically negligible. Moreover, note that the factors of $N$ obtained in this preparation step are bounded by $O(N^{1/3})$. Hence, we may compute their complete prime factorizations in negligible time by further applications of Corollary \ref{c1}. Let $N_k$ be the number which remains after these computations. Since no factor smaller or equal to $\lceil N_k^{1/3} \rceil$ has been found by applying Corollary \ref{c1}, $N_k$ does not have more than two prime factors. We may check easily if $N_k$ is a square number; if not, $N_k$ has to be a prime or a semiprime number. 

For the remainder of the proof, we may hence suppose that $N$ is either prime or $N=pq$ such that $p,q$ are distinct primes. We now run the following procedure.

\begin{algorithm}
	\emph{Input:} A prime or semiprime $N$, $\xi:=\lceil N^{1/9}/\log^{2/3} N\rceil$ and $\eta:=\lceil N^{2/9}/\log^{1/3} N\rceil$.
	
	\emph{Output:} The prime factorization of $N$.
	\begin{algorithmic}[1]
		\State Apply Corollary \ref{c1} with $\Delta=\lceil\left(N/(\eta+1)\right)^{1/2}\rceil$. If a factor is found, return and stop.
		\State Apply Theorem \ref{t3} with $\delta=\lceil N^{2/5}\rceil$. If a factor of $N$ is found or $N$ is proved to be prime, return and stop; otherwise, $\alpha\in\Z_N^*$ is found such that $\ord_N(\alpha)>\lceil N^{2/5}\rceil$.
		\State Apply Algorithm \ref{a2} with $\xi$, $\eta$ and the element $\alpha$ found in Step 2. If no factor is found, return that $N$ is prime.
	\end{algorithmic}
	\label{a3}
\end{algorithm}

Note that, for sufficiently large inputs $N$, we have $\ord_N(\alpha)>N^{2/5}\gg N^{1/2}/(\xi^{1/2}\eta)$ and the parameters $\xi$ and $\eta$ satisfy the assumptions of Algorithm \ref{a2}. In Step 1, we make sure that $p>\left(N/(\eta+1)\right)^{1/2}$ holds. As a result, Lemma \ref{l5} implies that the procedure above finds the factors of a semiprime number. Therefore, if we do not find any proper factors, $N$ must be prime. We conclude that the algorithm is correct and are left with the final task to analyze the runtime complexity.

The runtime of Step 1 is bounded by $O\left(\textsf{M}_{\integ}\left(\Delta^{1/2}\log N\right)\log N\right)\subseteq \widetilde{O}(N^{7/36})$, which is asymptotically neglible. The same is true for Step 2, which runs in time
\[
O\left(\frac{N^{1/5}\log^2 N}{\sqrt{\log\log N}} \right).
\]
We now show that our choice of $\xi$ and $\eta$ optimizes the runtime complexity of Step 3. Considering the three summands in Lemma \ref{l5}, one easily observes that the best possible choice for $\xi$ and $\eta$ is obtained by solving
\begin{equation}
\eta\log\eta=\frac{N^{1/2}}{\xi^{1/2}\eta}=\frac{N^{1/4}\xi^{3/4}\log\xi}{\eta^{1/2}}
\label{e4}
\end{equation}
under $O$-notation. Assuming that both $\log\eta$ and $\log \xi$ are in $O(\log N)$, we solve $\eta^2=N^{1/2}/(\xi^{1/2}\log N)$ and $\eta^{3/2}=N^{1/4}\xi^{3/4}$ for $\xi$ and $\eta$. It is easy to check that our choice of these parameters is a solution to these equations, and that the three terms in (\ref{e4}) are all equal to $N^{2/9}\log^{2/3}N$ under $O$-notation. As a result, the overall bit-complexity of Step 3 may be bounded by 
\[
O\left(\textsf{M}_{\integ}\left(N^{2/9}\log^{5/3}N\right)\log N\right),
\]
and the claim follows. 
\qed

\section*{Acknowledgment}
I want to thank an anonymous referee for the helpful and most valuable suggestions in the
report on the first version of this paper.


\end{document}